\documentclass[12pt]{article}
\usepackage{amssymb,latexsym,theorem}

\newcommand{\Diag}{\mathop{\mathrm{Diag}}\nolimits}
\newcommand{\Spec}{\mathop{\mathrm{Spec}}\nolimits}

\newcommand{\ind}{\mathop{\mathrm{ind}}\nolimits}
\newcommand{\hull}{\mathop{\mathrm{hull}}\nolimits}

\newcommand{\hgt}{\mathop{\mathrm{ht}}\nolimits}
\newcommand{\gr}{\mathop{\mathrm{gr}}\nolimits}

\newcommand{\SL}{\mathop{\mathit{SL}}\nolimits}

\newcommand{\Z}{\mathbb Z}

\newcommand{\C}{\mathcal C}

\newcommand{\qed}{\unskip\nobreak\hfill\hbox{ $\Box$}}

\newtheorem{Theorem}[subsection]{Theorem}
\newtheorem{Lemma}[subsection]{Lemma}

\newtheorem{Hypothesis}[subsection]{Hypothesis}
\theorembodyfont{\normalfont}
\newtheorem{Remark}[subsection]{Remark}

\newtheorem{Example}[subsection]{Example}

\begin{document}
\title{Finite good filtration dimension for modules over an algebra with
good filtration.}
\author{Wilberd van der Kallen}
\date{
\it \small Dedicated to Eric Friedlander on his $60^{\mathrm{th}}$ birthday
}
\maketitle
\sloppy
\begin{abstract}
Let $G$ be a connected
reductive linear algebraic group over a field $k$ of
characteristic $p>0$.
Let $p$ be large enough with respect to the root system. We show that
if a finitely generated commutative $k$-algebra $A$ with $G$-action
has good filtration,
then any noetherian $A$-module with compatible $G$-action has
finite good filtration dimension.
\end{abstract}

\section{Introduction}
Consider  a connected reductive linear algebraic  group
$G$ defined over a field $k$ of positive characteristic $p$.
We say that $G$ has the cohomological finite generation property (CFG)
if the following holds:
Let $A$ be a finitely generated commutative $k$-algebra on which $G$ acts
 rationally by $k$-algebra automorphisms. (So $G$ acts on $\Spec(A)$.)
Then the cohomology ring $H^*(G,A)$ is finitely generated
as a $k$-algebra. Here, as in
\cite[I.4]{Jantzen}, we use the cohomology introduced by
Hochschild, also known as `rational cohomology'.

In \cite{cohGrosshans} we have shown that $\SL_2$ over a field of positive
characteristic has property (CFG), and in \cite{reductive} we proved that
$\SL_3$ over a field of
characteristic two has property (CFG). We conjecture that every reductive
linear algebraic group has property (CFG). In this paper we show that this
is at least a good heuristic principle: We derive one of the consequences
of (CFG) for any
simply connected
semisimple linear algebraic  group $G$ that satisfies the
following

\begin{Hypothesis}\label{gfhypo}
Assume that for every fundamental weight $\varpi_i$ the symmetric algebra
$S^*(\nabla_G(\varpi_i))$ on the fundamental representation
$\nabla_G(\varpi_i)$
has a good filtration.
\end{Hypothesis}

Recall that this hypothesis is satisfied if $p\geq\max_i(
\dim(\nabla_G(\varpi_i)))$, by \cite[4.1(5) and 4.3(1)]{Andersen-Jantzen}.
This inequality is not necessary. For instance, $\SL_n$ satisfies the
hypothesis
for $n\leq5$, by \cite[Lemma 3.2]{cohGrosshans}. When $p=2$, the hypothesis
 does not hold
 for $\SL_n$ with
$n\geq6$, by \cite[3.3]{cohGrosshans}.

In the sequel let $G$ be a connected reductive linear algebraic group
over a field $k$ of
characteristic $p>0$
with simply connected
commutator subgroup for which hypothesis \ref{gfhypo} holds.
Let $A$ be a finitely generated commutative $k$-algebra on which $G$ acts
rationally by $k$-algebra automorphisms. Let $M$ be a noetherian $A$-module
on which $G$ acts compatibly. This means that the structure map $A\otimes M\to
M$ is a $G$-module map.
Our main result is

\begin{Theorem}\label{main}
If $A$ has good filtration, then $M$ has finite good filtration dimension and
each $H^i(G,M)$ is a noetherian $A^G$-module.
\end{Theorem}

When $A=k$ the theorem goes back to \cite{Friedlander-Parshall}
 and does not need hypothesis \ref{gfhypo}.
Unlike the proofs in \cite{cohGrosshans} and \cite{reductive},
the proof of our theorem
does not involve any cohomology of finite group schemes and is thus
independent of the work of Friedlander and Suslin \cite{Friedlander-Suslin}.
But without their work we would not have guessed the theorem.
For clarity we
will pull some material of \cite{cohGrosshans} free from
 finite group schemes.

\section{Recollections}
Some unexplained notations, terminology, properties, \ldots can be
found in \cite{Jantzen}.
We choose a Borel group $B^+=TU^+$  and the
opposite Borel group $B^-$. The roots of $B^+$ are positive.
If $\lambda\in X(T)$ is dominant, then $\ind_{B^-}^G(\lambda)$ is the
`dual Weyl module' or `costandard module'
$\nabla_G(\lambda)$ with highest weight $\lambda$. The formula
$\nabla_G(\lambda)=\ind_{B^-}^G(\lambda)$ just means that $\nabla_G(\lambda)$
is obtained from the Borel-Weil construction:
$\nabla_G(\lambda)$ equals $H^0(G/B^-,{\mathcal L})$ for a certain
 line bundle on the
flag variety $G/B^-$.
In a \emph {good filtration}
$0=V_{-1}\subseteq V_0 \subseteq V_1\ldots$ of a $G$-module
$V=\bigcup_iV_i$ the nonzero layers $V_i/V_{i-1}$ are of the form
$\nabla_G(\mu)$. As in \cite{vdkallen book}
we will actually also allow a layer to be a direct sum
of any number of copies of the same $\nabla_G(\mu)$,
cf. \cite[II.4.16 Remark 1]{Jantzen}.
This is much more convenient
when working with
infinite dimensional $G$-modules. It is shown in \cite{Friedlander can}
that a module of countable
dimension that has a good filtration in our sense also has a filtration
that is a good filtration
in the old sense. Note that
the module $M$ in our theorem has countable dimension.
It would do little harm to restrict to modules of countable dimension
throughout.

If $V$ is a $G$-module, and $m\geq-1$ is an integer so that
$H^{m+1}(G,\nabla_G(\mu)\otimes V)=0$ for all dominant $\mu$, then we say
as in \cite{Friedlander-Parshall}
that $V$ has \emph{good filtration dimension} at most $m$.
The case $m=0$ corresponds with $V$ having a good filtration.
And for $m\geq0$ it means that $V$ has a resolution
$$0\to V\to N_0 \to \cdots \to N_m\to 0$$ in which the $N_i$ have good
filtration, in our sense.
We say that $V$ has good filtration dimension precisely $m$,
notation $\dim_\nabla(V)=m$,
 if $m$ is
minimal so that $V$ has good filtration dimension at most $m$.
In that case $H^{i+1}(G,\nabla_G(\mu)\otimes V)=0$ for all dominant $\mu$
and all $i\geq m$. In particular $H^{i+1}(G,V)=0$ for $i\geq m$.
If there is no finite $m$ so that  $\dim_\nabla(V)=m$, then we put
$\dim_\nabla(V)=\infty$.

\subsection{Filtrations}
For simplicity assume also that $G$ is semisimple. (Until remark
\ref{sssuffice}.)
If $V$ is a $G$-module, and $\lambda$ is  a dominant weight,
then $V_{\leq\lambda}$ denotes the largest $G$-submodule all whose weights
$\mu$ satisfy $\mu\leq\lambda$ in the dominance partial order
\cite[II.1.5]{Jantzen}. For instance, $V_{\leq0}$ is the module of invariants
$V^G$.
Similarly $V_{<\lambda}$ denotes the largest $G$-submodule all whose weights
$\mu$ satisfy $\mu<\lambda$. As in \cite{vdkallen book}, we form the
$X(T)$-graded module
$$\gr_{X(T)} V=\bigoplus_{\lambda\in X(T)}V_{\leq\lambda}/V_{<\lambda}.$$
Each $V_{\leq\lambda}/V_{<\lambda}$, or $V_{\leq\lambda/<\lambda}$ for short,
 has a $B^+$-socle
$(V_{\leq\lambda/<\lambda})^U=V^U_\lambda$ of weight
$\lambda$. We always view $V^U$ as a $B^-$-module through restriction
(inflation)
along the
homomorphism $B^-\to T$.
Then $V_{\leq\lambda/<\lambda}$ embeds naturally in its `good filtration
hull' $\hull_\nabla(V_{\leq\lambda/<\lambda})
=\ind_{B^-}^GV^U_\lambda$.
This good filtration hull has the same $B^+$-socle
and by Polo it is the injective hull
in the category $\C_\lambda$
of $G$-modules $N$ that satisfy $N=N_{\leq\lambda}$.
Compare \cite[3.1.10]{vdkallen book}.

We convert the $X(T)$-graded module $\gr_{X(T)} V$ to a $\Z$-graded
module through an additive height function $\hgt:X(T)\to \Z$, defined by
$\hgt=2\sum_{\alpha>0}\alpha^\vee$,
the sum being over the positive roots.
(Our  $\hgt$ is twice the one used by Grosshans \cite{Grosshans contr},
because we prefer to get
even degrees rather than just integer degrees.)
The Grosshans graded module is now
$$\gr V=\bigoplus_{i\geq0}\gr_i V,$$
with $$\gr_i V=\bigoplus_{\hgt(\lambda)=i}V_{\leq\lambda/<\lambda}.$$
In other words, if one puts
$$V_{\leq i}:=\sum_{\hgt(\lambda)\leq i}V_{\leq\lambda},$$
then $\gr V$ is is the associated graded of the
filtration $V_{\leq0}\subseteq V_{\leq1}\cdots$.

Let us apply the above to our finitely generated commutative
$k$-algebra with $G$-action $A$.
The Grosshans graded algebra $\gr A$ embeds in a good filtration hull,
which Grosshans calls $R$, and
which we call $\hull_\nabla(\gr A)$,
$$\hull_\nabla(\gr A):=\ind_{B^-}^GA^U=
\bigoplus_i\bigoplus_{\hgt(\lambda)=i}
\hull_\nabla(A_{\leq\lambda}/A_{<\lambda}).$$
Grosshans \cite{Grosshans contr} shows that $A^U$, $\gr A$,
$\hull_\nabla(\gr A)$ are finitely
 generated $k$-algebras with $\hull_\nabla(\gr A)$ finite over $\gr A$.
Mathieu studied $\gr A$ and $\hull_\nabla(\gr A)$ earlier
in \cite{Mathieu G}.

\begin{Example}\label{multicone}
Consider the multicone \cite{Kempf Ramanathan}
$$k[G/U]:=\ind_U^Gk=\ind_{B^+}^G\ind_U^{B^+}k=\ind_{B^+}^Gk[T]=\bigoplus_{
\lambda \textrm{ \scriptsize dominant}}\nabla_G(\lambda).$$
It is its own Grosshans graded ring. Recall \cite{Kempf Ramanathan} that
it is  generated as a $k$-algebra by
the finite dimensional sum of the  $\nabla_G(\varpi_i)$, where $\varpi_i$
denotes the $i$th fundamental weight.
\end{Example}

\begin{Lemma}\label{extend}
Let $A$ have a good filtration, so that $\gr A=\hull_\nabla(\gr A)$.
Let $R=\oplus_i R_i$ be a graded algebra with $G$-action
such that $R_i=(R_i)_{\leq i}$. Then every $T$-equivariant
graded algebra homomorphism $R^U\to (\gr A)^U$ extends uniquely to a
$G$-equivariant graded algebra homomorphism $R\to (\gr A)$.
\end{Lemma}

\paragraph{Proof}
Use that $\hull_\nabla(\gr A)$  is an induced module.\qed

\subsection{A graded polynomial $G\times D$-algebra with good filtration}

We now extract a construction from \cite{cohGrosshans}. It is hidden in the
study
of a Hochschild-Serre spectral sequence which in the present situation would
correspond with the case where as normal subgroup one takes the trivial
subgroup!

As the algebra
$(\gr A)^U$ is finitely generated, it is also generated by finitely many
weight vectors. Consider one such weight vector $v$, say of weight $\lambda$.
Clearly $\lambda$ is dominant.
If $\lambda=0$, map a polynomial ring $P_v:=k[x]$ with trivial $G$-action to
 $\gr A$ by substituting $v$ for $x$. Also put $D_v:=1$.
Next assume $\lambda\neq0$. Let $\ell$ be the rank of $G$.
Define a $T$-action
on the $X(T)$-graded algebra
$$P=\bigotimes_{i=1}^{\ell}S^*(\nabla_G(\varpi_i))$$
by letting $T$ act on $\bigotimes_{i=1}^{\ell}S^{m_i}(\nabla_G(\varpi_i))$
through weight $\sum_im_i\varpi_i$. So now we have a $G\times T$-action on $P$.
Observe
that by our key hypothesis \ref{gfhypo} and the tensor product
property \cite[Ch.~G]{Jantzen}
the polynomial algebra $P$ has a good filtration for the $G$-action.
Let $D$ be the scheme
theoretic kernel of $ \lambda$.
So $D$ has character group
$X(D)=X(T)/\Z \lambda$ and $D=\Diag(X(T)/\Z \lambda)$ in the notations of
\cite[I.2.5]{Jantzen}.
The subalgebra $P^{1\times D}$ is a graded algebra with
good filtration such that its subalgebra
$P^{U\times D}$ contains a polynomial algebra on one generator $x$ of weight
$\lambda\times \lambda$. In fact, this polynomial subalgebra contains all the
weight vectors in $P^{U\times D}$ of weight $\mu\times\nu$ with $\hgt(\mu)\geq
\hgt(\nu)$. The other weight vectors in $P^{U\times D}$
also have weight of the form $\mu\times\nu$
with $\nu$ a multiple of $\lambda$. These other weight vectors
span an ideal in $P^{U\times D}$.
Now assume $A$ has a good filtration.
By lemma \ref{extend} one easily constructs
 a $G$-equivariant algebra homomorphism
$P^{1\times D}\to \gr A$ that maps $x$ to $v$.
Write it as $P_v^{1\times D_v}\to \gr A$, to stress the dependence on $v$.

As new $P$ we take the tensor product of the finitely many $P_v$
and as diagonalized group $D$ we take the direct product of the $D_v$.
Then we have a  graded algebra map $P^D\to \gr A$.
It is surjective because its image has good filtration
(\cite[Ch.~A]{Jantzen}) and contains $(\gr A)^U$.
The $G\times D$-algebra $P$ is an example of what we called  in
\cite{cohGrosshans}
a graded polynomial $G\times D$-algebra with good
filtration. We have proved

\begin{Lemma}
If $A$ has a good filtration, then there is a graded polynomial
$G\times D$-algebra $P$ with good
filtration and a graded $G$-equivariant surjection $P^D\to \gr A$.
\end{Lemma}

Now recall $M$ is a noetherian $A$-module on which $G$ acts compatibly, meaning
that the structure map $A\otimes M\to M$ is a map of $G$-modules.
Form the `semi-direct product ring' $A\ltimes M$ whose underlying $G$-module
is $A\oplus M$, with product given by
$(a_1,m_1)(a_2,m_2)=(a_1a_2,a_1m_2+a_2m_1)$. By Grosshans $\gr(A\ltimes M)$ is
a finitely generated algebra, so we get

\begin{Lemma}
$\gr M$ is a noetherian $\gr A$-module.
\end{Lemma}

This is of course very reminiscent of the proof of the lemma
\cite[Thm. 16.9]{Grosshans book} telling that
$M^G$ is a noetherian module over the finitely generated $k$-algebra
$A^G$. We will tacitly use
its counterpart for diagonalized actions, cf. \cite{Borsari-Santos},
\cite[I.2.11]{Jantzen}.


Taking things together we learn that if $A$ has a good filtration,
then $P\otimes_{P^D} \gr M$ is what we called in \cite{cohGrosshans}
a finite graded
$P$-module. Lemma \cite[Lemma 3.7]{cohGrosshans} then tells us

\begin{Lemma}
Let $A$ have good filtration.
Then $P\otimes_{P^D}\gr M$ has finite good filtration dimension and
each $H^i(G,P\otimes_{P^D}\gr M)$ is a noetherian $P^G$-module.
\end{Lemma}

Extend the $D$-action on $P$ to $P\otimes_{P^D}\gr M$ by using the trivial
action on the second factor. Then we have a $G\times D$-module structure
on $P\otimes_{P^D}\gr M$. As $D$ is diagonalized,
$P^D$ is a direct summand of $P$ as a $P^D$-module \cite[I.2.11]{Jantzen}
and $(P\otimes_{P^D}\gr
M)^{1\times D}=\gr M$ is a direct summand of the $G$-module
$P\otimes_{P^D}\gr M$.
It follows that $\gr M$ also has finite good filtration dimension and it
follows that each $H^i(G,P\otimes_{P^D}\gr M)^{1\times D}=H^i(G,\gr M)$
is a noetherian $P^{G\times D}$-module.
But the action of $P^{G\times D}$ on $\gr M$
factors through $(\gr A)^G$, so we see
that each $H^i(G,\gr M)$ is a noetherian $(\gr A)^G$-module.
And one always has $(\gr A)^G=(\gr_0 A)^G=A^G$.
We conclude

\begin{Lemma}\label{small}
Let $A$ have good filtration. Then $\gr M$ has finite good filtration
dimension and
each $H^i(G,\gr M)$ is a noetherian $A^G$-module.
\end{Lemma}

\section{Degrading}
We still have to get rid of the grading.
The filtration $M_{\leq 0}\subseteq M_{\leq 1}\cdots$ induces a filtration
of the Hochschild complex \cite[I.4.14]{Jantzen} whence a spectral sequence
$$E(M):E_1^{ij}=H^{i+j}(G,\gr_{-i}M)\Rightarrow H^{i+j}(G, M).$$
It lives in an unusual quadrant.

Assume that $A$ has good filtration. Then by Lemma \ref{small}
$E_1( M)$ is a finitely generated $A^G$-module.
 So the spectral sequence lives in only finitely many
bidegrees $(i,j)$. Thus there is
the same kind of convergence as one would have in a more common quadrant.

Choose $A^G$ as ring of operators to act on the spectral sequence
$E(M)$.  As $E_1(M)$
is a noetherian $A^G$-module, it easily follows (even without the spectral
sequence)
that $H^*(G,M)$ is a noetherian $A^G$-module.
To finish the proof of the theorem,
we note that $A\otimes k[G/U]$ is also a finitely generated
algebra with a good filtration
and that $M\otimes k[G/U]$ is a noetherian module over it.
So what we have just seen tells that $H^*(G,M\otimes k[G/U])$ is a noetherian
$(A\otimes k[G/U])^G$-module. In particular, there is an $m\geq-1$ so
that $H^{m+1}(G,M\otimes k[G/U])=0$.\qed

\begin{Remark}\label{sssuffice}
Somewhere along the way we made the simplifying assumption that $G$ is
semisimple. So for  the original $G$ we have now proved that $M$ has finite
good filtration dimension with respect to the commutator subgroup $H$ of $G$.
But that is the same as having  finite
good filtration dimension with respect to $G$. Also, the fact that $H^i(H,M)$
is a noetherian $A^H$-module implies that $H^i(G,M)$
is a noetherian $A^G$-module by taking invariants under the
diagonalizable center $Z(G)$.
\end{Remark}

\begin{Remark}\label{AG}
We did not prove that $M$ has a finite resolution by noetherian $A$-modules
with compatible $G$-action and good filtration. 
We do not know how to start. One may embed $M$ into
the $A$-module  $M\otimes k[G]$ with compatible $G$-action. 
It has good filtration, but it is not noetherian as an $A$-module.
\end{Remark}

\end{document}